\documentclass[11pt]{article}

\usepackage{amssymb,latexsym,amsmath}

\usepackage{graphicx}

\begin{document}

\newcommand{\bfi}{\bfseries\itshape}

\makeatletter

\@addtoreset{figure}{section}

\def\thefigure{\thesection.\@arabic\c@figure}

\def\fps@figure{h, t}

\@addtoreset{table}{bsection}

\def\thetable{\thesection.\@arabic\c@table}

\def\fps@table{h, t}

\@addtoreset{equation}{section}

\def\theequation{\thesubsection.\arabic{equation}}

\makeatother

\newtheorem{thm}{Theorem}[section]

\newtheorem{prop}[thm]{Proposition}

\newtheorem{lema}[thm]{Lemma}

\newtheorem{cor}[thm]{Corollary}

\newtheorem{defi}[thm]{Definition}

\newtheorem{rk}[thm]{Remark}

\newtheorem{exempl}{Example}[section]

\newenvironment{exemplu}{\begin{exempl}  \em}{\hfill $\square$

\end{exempl}}

\newcommand{\comment}[1]{\par\noindent{\raggedright\texttt{#1}

\par\marginpar{\textsc{Comment}}}}

\newcommand{\todo}[1]{\vspace{5 mm}\par \noindent \marginpar{\textsc{ToDo}}\framebox{\begin{minipage}[c]{0.95 \textwidth}

\tt #1 \end{minipage}}\vspace{5 mm}\par}

\newcommand{\ea}{\mbox{{\bf a}}}

\newcommand{\eu}{\mbox{{\bf u}}}

\newcommand{\ueu}{\underline{\eu}}

\newcommand{\ueo}{\overline{u}}

\newcommand{\oeu}{\overline{\eu}}

\newcommand{\ew}{\mbox{{\bf w}}}

\newcommand{\ef}{\mbox{{\bf f}}}

\newcommand{\eF}{\mbox{{\bf F}}}

\newcommand{\eC}{\mbox{{\bf C}}}

\newcommand{\en}{\mbox{{\bf n}}}

\newcommand{\eT}{\mbox{{\bf T}}}

\newcommand{\eL}{\mbox{{\bf L}}}

\newcommand{\eR}{\mbox{{\bf R}}}

\newcommand{\eV}{\mbox{{\bf V}}}

\newcommand{\eU}{\mbox{{\bf U}}}

\newcommand{\ev}{\mbox{{\bf v}}}

\newcommand{\eve}{\mbox{{\bf e}}}

\newcommand{\uev}{\underline{\ev}}

\newcommand{\eY}{\mbox{{\bf Y}}}

\newcommand{\eK}{\mbox{{\bf K}}}

\newcommand{\eP}{\mbox{{\bf P}}}

\newcommand{\eS}{\mbox{{\bf S}}}

\newcommand{\eJ}{\mbox{{\bf J}}}

\newcommand{\eB}{\mbox{{\bf B}}}

\newcommand{\eH}{\mbox{{\bf H}}}

\newcommand{\leb}{\mathcal{ L}^{n}}

\newcommand{\eI}{\mathcal{ I}}

\newcommand{\eE}{\mathcal{ E}}

\newcommand{\hen}{\mathcal{H}^{n-1}}

\newcommand{\eBV}{\mbox{{\bf BV}}}

\newcommand{\eA}{\mbox{{\bf A}}}

\newcommand{\eSBV}{\mbox{{\bf SBV}}}

\newcommand{\eBD}{\mbox{{\bf BD}}}

\newcommand{\eSBD}{\mbox{{\bf SBD}}}

\newcommand{\ecs}{\mbox{{\bf X}}}

\newcommand{\eg}{\mbox{{\bf g}}}

\newcommand{\paromega}{\partial \Omega}

\newcommand{\gau}{\Gamma_{u}}

\newcommand{\gaf}{\Gamma_{f}}

\newcommand{\sig}{{\bf \sigma}}

\newcommand{\gac}{\Gamma_{\mbox{{\bf c}}}}

\newcommand{\deu}{\dot{\eu}}

\newcommand{\dueu}{\underline{\deu}}

\newcommand{\dev}{\dot{\ev}}

\newcommand{\duev}{\underline{\dev}}

\newcommand{\weak}{\stackrel{w}{\approx}}

\newcommand{\mild}{\stackrel{m}{\approx}}

\newcommand{\strong}{\stackrel{s}{\approx}}

\newcommand{\weakdown}{\rightharpoondown}

\newcommand{\opg}{\stackrel{\mathfrak{g}}{\cdot}}

\newcommand{\opunu}{\stackrel{1}{\cdot}}
\newcommand{\opdoi}{\stackrel{2}{\cdot}}

\newcommand{\opn}{\stackrel{\mathfrak{n}}{\cdot}}
\newcommand{\opx}{\stackrel{x}{\cdot}}

\newcommand{\tr}{\ \mbox{tr}}

\newcommand{\Ad}{\ \mbox{Ad}}

\newcommand{\ad}{\ \mbox{ad}}

\renewcommand{\contentsname}{ }

\title{Dilatation structures with the Radon-Nikodym property}

\author{Marius Buliga \\
\\
Institute of Mathematics, Romanian Academy \\
P.O. BOX 1-764, RO 014700\\
Bucure\c sti, Romania\\
{\footnotesize Marius.Buliga@imar.ro}}

\date{This version:  25.06.2007}

\maketitle

\section*{Introduction}

The notion of a dilatation structure stemmed out from my efforts to
understand basic results in sub-Riemannian geometry, especially the last section
of the paper by Bella\"{\i}che \cite{bell} and the intrinsic point of view 
of Gromov \cite{gromo}. 

In these papers, as in other articles devoted to sub-Riemannian geometry,
fundamental results admiting an intrinsic formulation were proved using
differential geometry tools, which are in my opinion not intrinsic to
sub-Riemannian geometry. 

Therefore I tried to find a self-contained  frame in which sub-Riemannian
geometry would be a model, if we use the same manner of speaking as in the 
case of   hyperbolic 
geometry (with its self-contained collection of axioms) and the Poincar\'e 
disk as a model of hyperbolic geometry. 

An outcome of this effort are the notions of a dilatation structure and a pair 
of dilatation structures, one looking down to another. To the first notion are
dedicated the papers \cite{buligadil1}, \cite{buligacont} (the second paper
treating about a "linear" version of a generalized dilatation structure,
corresponding to Carnot groups or more general contractible groups). 

As it seems now, dilatation structures are a valuable notion by itself, with 
possible field of application strictly containing sub-Riemannian geometry, but 
also ultrametric spaces or contractible groups. A dilatation structure encodes 
the approximate self-similarity of a metric space and it induces non associative 
but approximately associative operations on the metric space, as well as 
a tangent bundle (in the metric sense) with group operations in each fiber 
(tangent space to a point). 

In this paper I explain what is a pair of dilatation structures, one looking 
down to another, see definition \ref{dlook}. Such a pair of dilatation
structures leads to the intrinsic definition of a distribution as a field of topological
filters, definition \ref{deftopdistri}.  

To any pair of dilatation structures there is an associated notion of
differentiability which generalizes the Pansu differentiability \cite{pansu}. 
This allows  the introduction of the Radon-Nikodym property for dilatation 
structures, which is the straightforward generalization of the Radon-Nikodym
property for Banach spaces. 

After an introducting section about length metric spaces and metric derivatives,
is proved in theorem \ref{fleng} that for a dilatation structure with the Radon-Nikodym property the 
length of absolutely continuous curves expresses as an integral of the norms 
of the tangents to the curve, as in Riemannian geometry.

Further it is shown  that Radon-Nikodym property transfers from any "upper" 
dilatation structure looking down  to a "lower" dilatation structure, theorem 
\ref{ttransfer}. Im my opinion this result explains intrinsically the fact that 
absolutely continuous curves in regular sub-Riemannian manifolds are derivable 
almost everywhere, as proved by Margulis, Mostow \cite{marmos}, Pansu
\cite{pansu} (for Carnot groups) or Vodopyanov \cite{vodopis}.  

The subject of application of these results for regular sub-Riemannian manifold 
will be left for a future paper, due to the unavoidable accumulation of 
technical estimates which are needed. 

\paragraph{Acknowledgements.}  I became accustomed with this subject of 
research during my period of work at the Department of Mathematics of the EPFL. 
During this time I had the opportunity to collaborate with some of the people
from the Institute of Mathematics of the Bern University, as well as with some 
invited guests at EPFL. For many valuable
discussions on the subject of sub-Riemannian geometry,  I would like to 
express my thanks and appreciation especially to Martin Reimann and Sergey 
Vodopyanov.

\tableofcontents

\section{Notations}

Let $\Gamma$ be  a topological separated commutative group  endowed with a continuous group morphism 
$$\nu : \Gamma \rightarrow (0,+\infty)$$ with $\displaystyle \inf \nu(\Gamma)  =  0$. Here $(0,+\infty)$ is 
taken as a group with multiplication. The neutral element of $\Gamma$ is denoted by $1$. We use the multiplicative notation for the operation in $\Gamma$. 

The morphism $\nu$ defines an invariant topological filter on $\Gamma$ (equivalently, an end). Indeed, 
this is the filter generated by the open sets $\displaystyle \nu^{-1}(0,a)$, $a>0$. From now on 
we shall name this topological filter (end) by "0" and we shall write $\varepsilon \in \Gamma \rightarrow 
0$ for $\nu(\varepsilon)\in (0,+\infty) \rightarrow 0$. 

The set $\displaystyle \Gamma_{1} = \nu^{-1}(0,1] $ is a semigroup. We note $\displaystyle 
\bar{\Gamma}_{1}= \Gamma_{1} \cup \left\{ 0\right\}$
On the set $\displaystyle 
\bar{\Gamma}= \Gamma \cup \left\{ 0\right\}$ we extend the operation on $\Gamma$ by adding the rules  
$00=0$ and $\varepsilon 0 = 0$ for any $\varepsilon \in \Gamma$. This is in agreement with the invariance 
of the end $0$ with respect to translations in $\Gamma$.

The space $(X,d)$ is a complete, locally compact metric space. For any $r>0$  
and any $x \in X$ we denote by $B(x,r)$ the open ball of center $x$ and radius 
$r$ in the metric space $X$.

By $\mathcal{O}(\varepsilon)$ 
we mean a positive function $f: \Gamma \rightarrow [0,+\infty)$  such that 
$\displaystyle \lim_{\varepsilon \rightarrow 0} f(\nu(\varepsilon)) \ = \ 0$.

\section{Length and metric derivatives}

For a detailed intrduction into the subject see for example \cite{amb}, chapter
1. 

\begin{defi}
The (upper) dilatation of a map $f: X \rightarrow Y$ between metric spaces,  in a point $u \in Y$ is 
$$ Lip(f)(u) = \limsup_{\varepsilon \rightarrow 0} \ \sup 
 \left\{ 
\frac{d_{Y}(f(v), f(w))}{d_{X}(v,w)} \ : \ v \not = w \ , \ v,w \in B(u,\varepsilon)
 \right\}$$
\end{defi}

In the particular case of a   derivable function $f: \mathbb{R} \rightarrow \mathbb{R}^{n}$ the upper dilatation is  $Lip(f)(t) \ = \ \mid \dot{f}(t) \mid$.  For any  Lipschitz function $f:X \rightarrow Y$ and for 
any $x \in X$ we have the obvious relation: 
$$Lip(f)(x) \ \leq \ Lip(f) \ .$$

A curve is a continuous function $f: [a,b] \rightarrow X$. The image of a curve is called path. Length measures paths. Therefore length does not depends on the reparametrisation of the path and it is additive with respect to concatenation of paths. 
 
In a metric space $(X,d)$ one can measure the length of curves in several ways.

\begin{defi}
The length 
of a  curve with $L^{1}$ dilatation $f: [a,b] \rightarrow X$ is 
$$L(f) = \int_{a}^{b} Lip(f)(t) \mbox{ d}t$$
\label{def1length}
\end{defi}

A different way to define a length of a curve is to consider its variation. 

\begin{defi}
The curve $f$ has bounded
variation if the quantity 
$$ Var(f) \ = \ \sup \left\{ \sum_{i=0}^{n} d(f(t_{i}),
f(t_{i+1})) \ \mbox{ : } a = t_{0} <  t_{1} < ... < t_{n} < t_{n+1} =  b
\right\}$$ 
(called variation of $f$) is finite. 
\label{def2length}
\end{defi}

There is a third, more basic way to introduce the length of a curve in a metric
space. 

\begin{defi}
The length of the path $A = f([a,b])$ is  the
one-dimensional Hausdorff measure of the path. The definition is the
following: 
$$l(A) \ = \ \lim_{\delta \rightarrow 0}  
 \inf \left\{ \sum_{i \in I} diam \ E_{i}  \mbox{ : } diam \ E_{i} 
< \delta \ , \ \ A \subset \bigcup_{i \in I} E_{i} \right\} $$
\label{def3lenght}
\end{defi}

The definitions are not equivalent. The variation $Var(f)$ of a curve $f$ and 
the length of a path L(f)  do not agree in general. Consider for example: 
$f: [-1,1] \rightarrow R^{2}$, $f(t) \ = \ (t, \ sign(t))$. 
We have $Var(f) \ = \ 4$ and $L(f([-1,1]) =  2$. Another example: 
the Cantor staircase function is continuous, but not Lipschitz. It has variation equal to 1 and 
length of the graph equal to 2. 

Nevertheless, for Lipschitz functions, the first  two definitions agree. For injective Lipschitz functions (i.e. for simple Lipschitz curves) the last two definitions agree.

\begin{thm}
For each Lipschitz curve $f: [a,b] \rightarrow X$, we have 
$L(f) \ = \ Var(f)$. 
\label{t411amb}
\end{thm}

\begin{thm} Suppose that $f: [a,b] \rightarrow X$ is a
 Lipschitz function and $A \ = \ f([a,b])$. Then $\displaystyle \mathcal{H}^{1}(A)  \leq  Var(f)$. 
 
 If $f$ is moreover injective then  $\displaystyle \mathcal{H}^{1}(A)  = Var(f)$. 
\label{t441amb}
\end{thm}

An important tool used in the proof  of the previous theorem is 
the geometrically obvious, but not straightforward to prove in this generality, 
Reparametrisation Theorem. 

\begin{thm}
Any path $A \subset X$ with a Lipschitz parametrisation admits a reparametrisation $f: [a,b] \rightarrow A$ such that $Lip(f)(t) = 1$ for almost any $t \in [a,b]$. 
\label{tp}
\end{thm}

We shall denote by $l_{d}$ the length functional, defined only on
Lipschitz curves, induced by the distance $d$. 
The length induces a new distance $d_{l}$, say on any Lipschitz connected component of the space $(X,d)$. The distance $d_{l}$ is given by: 
$$d_{l}(x,y) \  = \ \inf \ \left\{ l_{d}(f([a,b])) \mbox{ : } f: [a,b] \rightarrow X \ \mbox{ Lipschitz } , \right.$$
$$\left. \ f(a)=x \ , \ f(b) = y \right\}$$

We have therefore  two
operators $d \mapsto l_{d}$ and $l \mapsto d_{l}$. This leads to the introduction
of length metric spaces. 

\begin{defi}
A length metric space is a metric space $(X,d)$ such that $d  = d_{l}$. 
\label{dpath}
\end{defi}

From theorem \ref{t411amb} we deduce that Lipschitz curves in complete length metric spaces are absolutely continuous. Indeed, here is the definition of an absolutely continuous curve (definition 1.1.1, chapter 1,   \cite{amb}). 
 
 \begin{defi}
 Let $(X,d)$ be a complete metric space. A curve $c:(a,b)\rightarrow X$ is absolutely continuous if there exists $m\in L^{1}((a,b))$ such that for any $a<s\leq t<b$ we have 
 $$d(c(s),c(t)) \leq \int_{s}^{t} m(r) \mbox{ d}r   .$$
 Such a function $m$ is called an upper gradient of the curve $c$. 
 \label{defac}
 \end{defi}
 
 According to theorem \ref{t411amb}, for a Lipschitz curve $c:[a,b]\rightarrow X$ in a complete lenght metric space such a function 
 $m\in L^{1}((a,b))$  is the upper dilatation  $Lip(c)$. 
More can be said about the expression of the upper dilatation. We need first to introduce the notion of 
metric derivative of a Lipschitz curve. 

\begin{defi}
A curve $c:(a,b)\rightarrow X$ is metrically derivable in $t\in(a,b)$ if the limit 
$$md(c)(t) = \lim_{s\rightarrow t} \frac{d(c(s),c(t))}{\mid s-t \mid}$$
exists and it is finite. In this case $md(c)(t)$ is called the metric derivative of $c$ in $t$. 
\label{defmd}
\end{defi}

For the proof of the following theorem see \cite{amb}, theorem 1.1.2, chapter 1. 

\begin{thm}
Let $(X,d)$ be a complete metric space and $c:(a,b)\rightarrow X$ be an absolutely continuous curve. 
Then $c$ is metrically  derivable for $\mathcal{L}^{1}$-a.e. $t\in(a,b)$. Moreover the function $md(c)$ belongs to $L^{1}((a,b))$ and it is minimal in the following sense: $md(c)(t)\leq m(t)$  for  $\mathcal{L}^{1}$-a.e. $t\in(a,b)$, for each upper gradient $m$ of the curve $c$. 
\label{tupper}
\end{thm}

\section{The Radon-Nikodym property}
\label{radon}
 
 \begin{defi} 
A dilatation structure $(X,d, \delta)$ has the Radon-Nikodym property if  any  
Lipschitz curve $\displaystyle c : [a,b] \rightarrow (X,d)$ is derivable 
almost everywhere. 
 \end{defi}
 
 \begin{exemplu}
For  $\displaystyle (X,d)  =  ( \mathbb{V}, d)$, a real, finite dimensional,
normed vector space, with distance $d$ induced by the norm, the (usual) 
 dilatations $\displaystyle \delta^{x}_{\varepsilon}$ are given by:  
$$ \delta_{\varepsilon}^{x} y \ = \ x + \varepsilon (y-x) $$
Dilatations are defined everywhere. The group $\Gamma$ is $(0,+\infty)$ and 
the function $\nu$ is the identity. 

There are few things to check (see the appendix):  axioms 0,1,2 are obviously 
true. For axiom A3, remark that for any $\varepsilon > 0$, $x,u,v \in X$ we 
have: 
$$\frac{1}{\varepsilon} d(\delta^{x}_{\varepsilon} u , 
\delta^{x}_{\varepsilon} v ) \ = \ d(u,v) \ , $$
therefore for any $x \in X$ we have $\displaystyle d^{x} = d$. 

Finally, let us check the axiom A4. For any $\varepsilon > 0$ and $x,u,v \in X$ we have
$$\delta_{\varepsilon^{-1}}^{\delta_{\varepsilon}^{x} u} \delta_{\varepsilon}^{x} v \ = \ 
x + \varepsilon  (u-x) + \frac{1}{\varepsilon} \left( x+ \varepsilon(v-x) - x - \varepsilon(u-x) \right) \ = \ $$
$$ = \ x + \varepsilon  (u-x) + v - u$$ 
therefore this quantity converges to 
$$x + v - u \ = \ x + (v - x) - (u - x)$$
as $\varepsilon \rightarrow 0$. The axiom A4 is verified. 

This dilatation structure has the Radon-Nikodym property. 
\label{401}
\end{exemplu}

\begin{exemplu}
Because dilatation structures are defined by local requirements, we can easily define dilatation 
structures on riemannian manifolds, using particular atlases of the manifold and the riemannian 
distance (infimum of length of curves joining two points). Note that any finite dimensional manifold can be endowed with a riemannian metric. This class of examples covers all dilatation structures used in 
differential geometry. The axiom A4 gives an operation of addition of vectors 
in the tangent space  (compare with Bella\"{\i}che \cite{bell} last section).  
\label{411riemann}
\end{exemplu}

\begin{exemplu} Take $\displaystyle X = \mathbb{R}^{2}$ with the euclidean 
distance $\displaystyle d$. For any $z \in \mathbb{C}$ of the 
form $z= 1+ i \theta$ we define dilatations 
$$\delta_{\varepsilon} x = \varepsilon^{z} x  \ .$$
It is easy to check that $\displaystyle (\mathbb{R}^{2},d, \delta)$ 
is a dilatation structure, with  dilatations 
$$\delta^{x}_{\varepsilon} y = x + \delta_{\varepsilon} (y-x)  \quad .$$

Two such dilatation structures (constructed with the help of complex numbers 
$1+ i \theta$ and $1+ i \theta'$) are equivalent if and only if $\theta = \theta'$.  

There are two other interesting  properties of these dilatation structures. 
The first is that if $\theta \not = 0$ then there are no non trivial 
Lipschitz curves in $X$ which are differentiable almost everywhere. It means
that such dilatation structure does not have the Radon-Nikodym property. 

The second property is that any holomorphic and Lipschitz function from $X$ to $X$ (holomorphic in the 
usual sense on $X = \mathbb{R}^{2} = \mathbb{C}$) is differentiable almost everywhere, but there are 
Lipschitz functions from $X$ to $X$ which are not differentiable almost everywhere (suffices to take a 
$\displaystyle \mathcal{C}^{\infty}$ function from  $\displaystyle \mathbb{R}^{2}$ to $\displaystyle \mathbb{R}^{2}$ which is not holomorphic). 
\label{exemp1}
\end{exemplu}

 The Radon-Nikodym property can be stated in two equivalent ways. 
 
 \begin{prop}
 Let $(X,d,\delta)$ be a  dilatation structure. Then the following are
 equivalent: 
 \begin{enumerate}
 \item[(a)] $(X,d,\delta)$ has the Radon-Nikodym property; 
 \item[(b)] any  
Lipschitz curve $\displaystyle c' : [a',b'] \rightarrow (X,d)$ admits a
reparametrization  $\displaystyle c : [a, b] \rightarrow (X,d)$ such that 
  for almost every $t \in [a,b]$ there is $\dot{c}(t) \in U(c(t))$ such that 
$$\frac{1}{\varepsilon} d(c(t+\varepsilon) , \delta_{\varepsilon}^{c(t)} \dot{c}(t)) \rightarrow 0 $$
$$\frac{1}{\varepsilon} d(c(t-\varepsilon) , \delta_{\varepsilon}^{c(t)} 
inv^{c(t)}( \dot{c}(t))) \rightarrow 0  \quad ; $$
\item[(c)] any  
Lipschitz curve $\displaystyle c' : [a',b'] \rightarrow (X,d)$ admits a
reparametrization  $\displaystyle c : [a, b] \rightarrow (X,d)$ such that 
 for almost every $t \in [a,b]$ there is a conical group morphism 
 $$\dot{c}(t) : \mathbb{R} \rightarrow T_{c(t)} X$$ 
 such that for any $a \in \mathbb{R}$ we have 
 $$\frac{1}{\varepsilon} d(c(t+\varepsilon a) , \delta_{\varepsilon}^{c(t)} \dot{c}(t)(a)) \rightarrow 0  . $$
\end{enumerate}
\end{prop}

\paragraph{Proof.}   
 It is  straightforward that a conical group morphism $f: \mathbb{R} \rightarrow (N,\delta)$ is defined by its value $f(1)\in N$. Indeed, for any $a>0$ we have $\displaystyle f(a) = \delta_{a} f(1)$ and for any 
 $a<0$ we have $\displaystyle f(a) = \delta_{a} f(1)^{-1}$. From the morphism property we also 
 deduce that  
 $$\delta v = \left\{ \delta_{a} v \mbox{ : } a>0 , v=f(1) \mbox{ or } v=f(1)^{-1} \right\}$$
 is a one parameter group and that for all $\alpha, \beta >0$ we have 
 $$\delta_{\alpha+\beta} u = \delta_{\alpha}u  \delta_{\beta}u  \quad  \square $$

\begin{defi}
In a conical group $N$   we shall denote by $D(N)$ the set of all $u\in N$ with the property that 
$\displaystyle \varepsilon \in ((0,\infty),+) \mapsto \delta_{\varepsilon} u \in N$ is a morphism of semigroups .

 $D(N)$ is always non empty, because it contains the neutral element of $N$. $D(N)$ is also a cone,  with dilatations $\displaystyle   \delta_{\varepsilon}$, and a closed set. 
\label{defdisn}
\end{defi}
  
  We shall always  identify a conical  group morphism $f:\mathbb{R} \rightarrow
  N$ with its value $f(1) \in D(N)$.

\subsection{Length formula from Radon-Nikodym property}

 \begin{thm}
Let $(X,d,\delta)$ be a dilatation structure with the Radon-Nikodym property, over a complete length metric space $(X,d)$. Then   for any Lipschitz curve $c:[a,b]\rightarrow X$ the  length of $\gamma= c([a,b])$ is 
$$L(\gamma) = \int_{a}^{b} d^{c(t)}(c(t),\dot{c}(t)) \mbox{ d}t  . $$
\label{fleng}
\end{thm}

\paragraph{Proof.}
The upper dilatation of $c$ in $t$ is 
$$ Lip(c)(t) = \limsup_{\varepsilon \rightarrow 0} \ \sup 
 \left\{ 
\frac{d(c(v), c(w))}{\mid v-w\mid} \ : \ v \not = w \ , \ \mid v-t\mid, \mid w-t\mid <  \varepsilon \right\}  . $$
From theorem \ref{tupper}Ê we deduce that for almost every $t\in(a,b)$ we have 
$$Lip(c)(t) = \lim_{s\rightarrow t} \frac{d(c(s),c(t))}{\mid s-t \mid} . $$
 
If the dilatation structure has the Radon-Nikodym property then for almost every $t \in [a,b]$ there is 
$\displaystyle \dot{c}(t) \in D(T_{c(t)} X)$ such that 
$$\frac{1}{\varepsilon} d(c(t+\varepsilon) , \delta_{\varepsilon}^{c(t)} \dot{c}(t)) \rightarrow 0 .$$
Therefore for almost every $t \in [a,b]$ we have
$$Lip(c)(t) = \lim_{\varepsilon\rightarrow 0} \frac{1}{\varepsilon} d(c(t+\varepsilon),c(t)) = d^{c(t)}(c(t),\dot{c}(t))  . $$
The formula for length follows from here. \quad $\square$

\subsection{A dilatation structure looking down to another}

Consider two dilatation structures $\displaystyle \mathcal{A} = (X,d_{A}, \delta)$ and 
$\displaystyle \mathcal{B} =  (X,d_{B}, \bar{\delta})$. We explain here in which sense  
$\mathcal{A}$ looks down at $\mathcal{B}$.

\begin{defi}
Given dilatation structures      $\displaystyle \mathcal{A} = (X,d_{A}, \delta)$ and 
$\displaystyle \mathcal{B} =  (X,d_{B}, \bar{\delta})$, we write that $\mathcal{A} \geq \mathcal{B}$ if the following conditions are fulfilled:  
\begin{enumerate} 
\item[(a)] the identity $\displaystyle id \ : (X,d_{A}) \rightarrow (X,d_{B})$ is 1-Lipschitz, 
\item[(b)] the identity $\displaystyle id \ : (X,d_{A}) \rightarrow (X,d_{B})$ is derivable 
everywhere and for any point $x \in X$ the derivative $D \ id(x)$ is a projector, 
\item[(c)] for any $x \in X$, any continuous curve 
$\varepsilon \in [0,1) \mapsto z(\varepsilon) \in X$, such that $\displaystyle 
d^{x}_{A}(z(0), x) \leq 3/2$, if  
$$\lim_{\varepsilon \rightarrow 0} \left( d^{x}_{A}(x, z(\varepsilon)) - 
\frac{1}{\varepsilon} \, d_{B}^{x}(x, \delta^{x}_{\varepsilon} z(\varepsilon))
\right) \ = \ 0$$  
then $\displaystyle \lim_{\varepsilon \rightarrow 0} 
d^{x}_{A}(Q_{\varepsilon}^{x} z(\varepsilon),
z(\varepsilon)) \ = 0$, where $\displaystyle Q_{\varepsilon}^{x} = \ \bar{\delta}^{x}_{\varepsilon^{-1}} \delta_{\varepsilon}^{x}$.  
\end{enumerate}
\label{dlook}
\end{defi}

We explain in more detail the meaning of this definition. Condition (a) says that for any $x,y\in X$ we have $\displaystyle d_{B}(x,y)\leq d_{A}(x,y)$. Condition (b) can be understood by using definition 
\ref{defdif}: for any $x\in X$ there  exists a 
 function $\displaystyle D \ id(x)$  defined on a neighbourhood of $x$ with values in  a neighbourhood 
 of $f(x)$ such that 
\begin{equation}
\lim_{\varepsilon \rightarrow 0} \sup \left\{  \frac{1}{\varepsilon} \overline{d} \left( \delta^{x}_{\varepsilon} u ,  \overline{\delta}^{x}_{\varepsilon} D \ id(x) (u) \right) \mbox{ : } d(x,u) \leq \varepsilon \right\}Ê  = 0 .
\label{dcondb}
\end{equation}
From here we deduce that  for any $x$ and $u$ such that $\displaystyle d_{B}(x,u)$ is sufficiently small 
$$ \lim_{\varepsilon \rightarrow 0} \overline{\delta}^{x}_{\varepsilon^{-1}} \delta^{x}_{\varepsilon} (u) = 
D \ id(x) (u)  $$
and the limit is uniform with respect to $u$. 

The second part of the condition (b) states that 
$$D \ id(x) D \ id(x) = D \ id(x) .$$
In order to understand the condition (c) we need to introduce the following 
topological version of a distribution. 

\begin{defi}
We denote by $TopD(x)$ the topological filter 
induced by the relatively open neighbourhoods of $x$ in the closed ball 
$\displaystyle \left\{ z \in X \mbox{ : } d^{x}_{A}(x,z) \leq 2 \right\}$, 
 given by 
$$F(x,\varepsilon, \lambda) \ = \ \left\{ z \in X \mbox{ : } d^{x}_{A}(x,z) \leq
2 \, \, , \, d^{x}_{A}(x, z) - 
\frac{1}{\varepsilon} \, d_{B}^{x}(x, \delta^{x}_{\varepsilon} z)
\leq \lambda \right\} \quad .$$
This filter is called the topological distribution associated with the pair 
of dilatation structures $\displaystyle \mathcal{A} = (X,d_{A}, \delta)$ and 
$\displaystyle \mathcal{B} =  (X,d_{B}, \bar{\delta})$, 
such that  $\mathcal{A} \geq \mathcal{B}$. 
\label{deftopdistri}
\end{defi}

With this notation we may rewrite the condition (c) definition \ref{dlook} like
this: let $z(\varepsilon)$ be a continuous curve such that (in the sense of
topological filters) 
$$ \lim_{\varepsilon \rightarrow 0} z(\varepsilon) \ \in \ TopD(x) \quad .$$ 
Then $\displaystyle \lim_{\varepsilon \rightarrow 0} 
d^{x}_{A}(Q_{\varepsilon}^{x} z(\varepsilon),
z(\varepsilon)) \ = 0$, where $\displaystyle Q_{\varepsilon}^{x} = \ \bar{\delta}^{x}_{\varepsilon^{-1}} \delta_{\varepsilon}^{x}$.  
This means that the "size of the vertical part" of $z(\varepsilon)$, which 
is $\displaystyle d^{x}_{A}(Q_{\varepsilon}^{x} z(\varepsilon),
z(\varepsilon))$, becomes arbitrarily small as $\varepsilon \rightarrow 0$.

\subsection{Transfer of Radon-Nikodym property}

Suppose that $\displaystyle (X,d_{A})$ and $\displaystyle (X,d_{B})$ are complete, locally compact, 
length metric spaces and that we have two dilatation structures      $\displaystyle \mathcal{A} = (X,d_{A}, \delta)$ and 
$\displaystyle \mathcal{B} =  (X,d_{B}, \bar{\delta})$, such that $\mathcal{A} \geq \mathcal{B}$. 

A sufficient condition to have (a) in definition \ref{dlook} is the following (true in the case of sub-Riemannian manifolds): 
\begin{enumerate}
\item[(a')] for any Lipschitz curve $c$, if $l_{A}(c) < + \infty$ then $l_{B}(c) = l_{A}(c)$. Here 
$l_{A}$ and $l_{B}$ denote the length functional associated to distance $d_{A}$, distance $d_{B}$ respectively. 
\end{enumerate}

We prove here the following result concerning the transfer of Radon-Nikodym property. 

\begin{thm}
Let $\displaystyle (X,d_{A})$ and $\displaystyle (X,d_{B})$ be complete, locally compact, 
length metric spaces. Suppose  that we have two dilatation structures      $\displaystyle \mathcal{A} = (X,d_{A}, \delta)$ and 
$\displaystyle \mathcal{B} =  (X,d_{B}, \bar{\delta})$, such that $\mathcal{A} \geq \mathcal{B}$. 
Under the assumptions (a') and  (b), (c), (d) from definition \ref{dlook}, if the dilatation structure 
$\displaystyle \mathcal{B} =  (X,d_{B}, \bar{\delta})$ has the Radon-Nikodym property, then the dilatation structure   $\displaystyle \mathcal{A} = (X,d_{A}, \delta)$  has the 
Radon-Nikodym property.
\label{ttransfer}
\end{thm}

\paragraph{Proof.}
Let $c: [0,1] \rightarrow (X,d_{A})$ be a Lipschitz curve. Because of hypothesis (a) it follows that 
$c: [0,1] \rightarrow (X,d_{B})$ is also Lipschitz. Moreover, we can reparametrize the curve $c$ with the 
$d_{A}$ lenght and so we can suppose that $c$ is $d_{A}$ 1-Lipschitz. Therefore we can suppose that $c$ is $d_{B}$ 1-Lipschitz. 

The dilatation structure $\displaystyle \mathcal{B} =  (X,d_{B}, \bar{\delta})$ has the 
Radon-Nikodym property. Then for  almost any $t \in [0,1]$ there is $\dot{c}(t)$ such that 
\begin{equation}
\frac{1}{\varepsilon} d_{B}(c(t+\varepsilon), \bar{\delta}^{x}_{\varepsilon} \dot{c}(t)) \rightarrow 0 \ \ \mbox{ as } \varepsilon \rightarrow 0
\label{eq1}
\end{equation}

\begin{equation}
\frac{1}{\varepsilon} d_{B}(c(t-\varepsilon), \bar{\delta}^{x}_{\varepsilon} \dot{c}(t)^{-1}) \rightarrow 0 \ \ \mbox{ as } \varepsilon \rightarrow 0
\label{eq1'}
\end{equation}

Further we shall give only the half of the proof, namely we shall use only relation \eqref{eq1}. To get a complete proof, one has to repeat the reasoning starting from \eqref{eq1'}.

Because $c$ is $d_{A}$  1- Lipschitz, it follows that 
$$d_{A}^{c(t)}(\delta_{\varepsilon^{-1}}^{c(t)} c(t+\varepsilon), \dot{c}(t)) \leq 2$$
 for any $\varepsilon < \varepsilon(t) \in (0,+\infty)$. From the local 
 compactness  with respect to $\displaystyle d_{A}^{c(t)}$ we find that 
for any $t \in [0,1]$ there is a sequence 
$\displaystyle (\varepsilon_{h})_{h} \subset (0,+ \infty)$, converging to $0$ 
as $h \rightarrow \infty$, and $u(t) \in X$ such that: 
$$\lim_{h \rightarrow \infty} \delta_{\varepsilon_{h}^{-1}}^{c(t)}   c(t+\varepsilon_{h}) \ = \ u(t)$$
Use equation \eqref{eq1} to get that 
$$\lim_{h \rightarrow \infty} \bar{\delta}_{\varepsilon_{h}^{-1}}^{c(t)}   c(t+\varepsilon_{h}) \ = \ \dot{c}(t)$$
Re-write this latter equation as: 
$$\lim_{h \rightarrow \infty} \bar{\delta}_{\varepsilon_{h}^{-1}}^{c(t)}  \delta_{\varepsilon_{h}}^{c(t)} \delta_{\varepsilon_{h}^{-1}}^{c(t)}  c(t+\varepsilon_{h}) \ = \ \dot{c}(t)$$
and use the first part of  hypothesis (b) to get 
$$D \ id (c(t)) u(t) \ = \ \dot{c}(t)$$

But according to the second part of the hypothesis (b) the operator $D \ id(c(t))$ is a projector, hence
$$D \ id (c(t)) \dot{c}(t) \ = \ \dot{c}(t)$$
Because of the fact that the derivative commutes with dilatations we get the important fact that for any 
$\varepsilon> 0$ 
\begin{equation}
\delta_{\varepsilon}^{c(t)} \dot{c}(t) \ = \ \bar{\delta}_{\varepsilon}^{c(t)} \dot{c}(t)
\label{eq2}
\end{equation}

We wish to prove
\begin{equation}
\frac{1}{\varepsilon} d_{A}(\delta_{\varepsilon}^{c(t)} \bar{\delta}_{\varepsilon^{-1}}^{c(t)} c(t+\varepsilon), 
c(t+\varepsilon)) \rightarrow 0 \ \ \mbox{ as } \varepsilon \rightarrow 0
\label{eq3}
\end{equation}
Suppose that \eqref{eq3} is true. Then we would have 
$$d_{A}^{c(t)} (\bar{\delta}_{\varepsilon^{-1}}^{c(t)} c(t+\varepsilon), \delta_{\varepsilon^{-1}}^{c(t)} c(t+ \varepsilon)) \rightarrow 0 \ \ \mbox{ as } \varepsilon \rightarrow 0$$

But relations  \eqref{eq1} and \eqref{eq2} imply  that 
$$d_{A}^{c(t)}(\bar{\delta}_{\varepsilon^{-1}}^{c(t)} c(t+\varepsilon) , \dot{c}(t) ) \rightarrow 0 \ \ \mbox{ as } \varepsilon \rightarrow 0$$
therefore we would finally get 
$$d_{A}^{c(t)}(\delta_{\varepsilon^{-1}}^{c(t)} c(t+\varepsilon) , \dot{c}(t) ) \rightarrow 0 \ \ \mbox{ as } \varepsilon \rightarrow 0$$
which is what we want  to prove: that the curve $c$ is derivable in $t$ with respect to the dilatation structure $\mathcal{A}$.  

Let us prove the relation  \eqref{eq3}. According to hypothesis (a') we have: 
$$0 \leq \frac{1}{\varepsilon} d_{A}(c(t+\varepsilon), c(t)) - \frac{1}{\varepsilon} d_{B} (c(t+\varepsilon), c(t)) \ \leq \  \frac{1}{\varepsilon} \int_{t}^{t+\varepsilon} \mid \dot{c}(\tau) \mid_{B} \mbox{ d}\tau \ - \ \frac{1}{\varepsilon} d_{B} (c(t+\varepsilon), c(t))$$
where the quantity 
$$\mid \dot{c}(s) \mid_{B} \ = \ \lim_{\varepsilon \rightarrow 0} \frac{d_{B}((c(s+\varepsilon), c(s))}{\varepsilon} \ = \ 
d_{B}^{c(s)} (c(s), \dot{c}(s))$$
exists for almost every $s \in [0,1]$, according to theorem \ref{tupper}. 

We obtain therefore the relation: 
\begin{equation}
\frac{1}{\varepsilon} d_{A}(c(t+\varepsilon), c(t)) - \frac{1}{\varepsilon} d_{B} (c(t+\varepsilon), c(t)) \ \rightarrow 0 \ \ \mbox{ as } \varepsilon \rightarrow 0
\label{eq4}
\end{equation}

Here is the moment to use the last hypothesis (d). Indeed, the  relation \eqref{eq4} implies that 
\begin{equation}
 d_{A}^{c(t)}(c(t), \delta_{\varepsilon^{-1}}^{c(t)}c(t+\varepsilon)) - \frac{1}{\varepsilon} d_{B}^{c(t)}(c(t+\varepsilon), c(t)) \ \rightarrow 0 \ \ \mbox{ as } \varepsilon \rightarrow 0
\label{eq4'}
\end{equation}
Denote by 
$\displaystyle z(t, \varepsilon) \ = \ 
\delta_{\varepsilon^{-1}}^{c(t)}c(t+\varepsilon)$. The relation \eqref{eq4'} 
becomes: 
\begin{equation}
 d_{A}^{c(t)}(c(t), z(t, \varepsilon)) - 
 \frac{1}{\varepsilon} d_{B}^{c(t)}(c(t), \delta_{\varepsilon}^{x}z(t,
 \varepsilon)) \ \rightarrow 0 \ \ \mbox{ as } \varepsilon \rightarrow 0
\label{eq4''}
\end{equation}
We also have 
$$d_{A}^{c(t)} (c(t), z(t, \varepsilon)) \ = \ 
\frac{1}{\varepsilon} d_{A}^{c(t)}(c(t), c(t+\varepsilon)) \leq 2$$
for $\varepsilon$ sufficiently small, because we supposed that $c$ was 
reparametrized with the length. 
Therefore, with the notations from definition \ref{deftopdistri} and the 
paragraph following it,  we have 
$$\lim_{\varepsilon \rightarrow 0} z(t, \varepsilon) \ \in  \ TopD(c(t))$$  
 From the hypothesis (d) we deduce that 
$$\lim_{\varepsilon \rightarrow 0} d_{A}^{c(t)}(z(t,\varepsilon), 
Q_{\varepsilon}^{c(t)} z(t,\varepsilon)) \ = \ 0$$
Let us see what this means: 
$$\lim_{\varepsilon \rightarrow 0} d_{A}^{c(t)} (\bar{\delta}^{c(t)}_{\varepsilon} c(t+\varepsilon), 
\delta_{\varepsilon}^{c(t)} c(t+\varepsilon)) \ = \ 0$$
This relation is equivalent with \eqref{eq3}, so the proof is done. \quad $\blacksquare$

\section{Appendix: Dilatation structures}

For the sake of completeness we list in this appendix the definition and
properties of a dilatation structure, according to \cite{buligadil1},
\cite{buligacont}.

\subsection{The axioms of a dilatation structure}

The  axioms of  a dilatation structure $(X,d,\delta)$ are listed further. 
The first axiom is merely a preparation for the next axioms. That is why we 
counted it as axiom 0.

\begin{enumerate}
\item[{\bf A0.}] The dilatations $$ \delta_{\varepsilon}^{x}: U(x) 
\rightarrow V_{\varepsilon}(x)$$ are defined for any 
$\displaystyle \varepsilon \in \Gamma, \nu(\varepsilon)\leq 1$. The sets 
$\displaystyle U(x), V_{\varepsilon}(x)$ are open neighbourhoods of $x$.  
All dilatations are homeomorphisms (invertible, continuous, with 
continuous inverse). 

We suppose  that there is a number  $1<A$ such that for any $x \in X$ we have 
$$\bar{B}_{d}(x,A) \subset U(x)  \ .$$
 We suppose that for all $\varepsilon \in \Gamma$, $\nu(\varepsilon) \in 
(0,1)$, we have 
$$ B_{d}(x,\nu(\varepsilon)) \subset \delta_{\varepsilon}^{x} B_{d}(x,A) 
\subset V_{\varepsilon}(x) \subset U(x) \ .$$

 There is a number $B \in (1,A)$ such that  for 
 any $\nu(\varepsilon) \in (1,+\infty)$ the associated dilatation  
$$\delta^{x}_{\varepsilon} : W_{\varepsilon}(x) \rightarrow B_{d}(x,B) \ , $$
is injective, invertible on the image. We shall suppose that 
$\displaystyle  W_{\varepsilon}(x)$ is a open neighbourhood of $x$, 
$$V_{\varepsilon^{-1}}(x) \subset W_{\varepsilon}(x) $$
and that for all $\displaystyle \varepsilon \in \Gamma_{1}$ and 
$\displaystyle u \in U(x)$ we have
$$\delta_{\varepsilon^{-1}}^{x} \ \delta^{x}_{\varepsilon} u \ = \ u \ .$$
\end{enumerate}

We have therefore  the following string of inclusions, for any $\varepsilon \in \Gamma$, $\nu(\varepsilon) \leq 1$, and any $x \in X$:
$$ B_{d}(x,\nu(\varepsilon)) \subset \delta^{x}_{\varepsilon}  B_{d}(x, A) 
\subset V_{\varepsilon}(x) \subset 
W_{\varepsilon^{-1}}(x) \subset \delta_{\varepsilon}^{x}  B_{d}(x, B) \quad . $$

A further technical condition on the sets  $\displaystyle V_{\varepsilon}(x)$ and $\displaystyle W_{\varepsilon}(x)$  will be given just before the axiom A4. (This condition will be counted as part of 
axiom A0.)

\begin{enumerate}
\item[{\bf A1.}]  We  have 
$\displaystyle  \delta^{x}_{\varepsilon} x = x $ for any point $x$. We also have $\displaystyle \delta^{x}_{1} = id$ for any $x \in X$.

Let us define the topological space
$$ dom \, \delta = \left\{ (\varepsilon, x, y) \in \Gamma \times X \times X 
\mbox{ : } \quad \mbox{ if } \nu(\varepsilon) \leq 1 \mbox{ then } y \in U(x) \,
\, , 
\right.$$ 
$$\left. \mbox{  else } y \in W_{\varepsilon}(x) \right\} $$ 
with the topology inherited from the product topology on 
$\Gamma \times X \times X$. Consider also $\displaystyle Cl(dom \, \delta)$, 
the closure of $dom \, \delta$ in $\displaystyle \bar{\Gamma} \times X \times X$ with product topology. 
The function $\displaystyle \delta : dom \, \delta \rightarrow  X$ defined by 
$\displaystyle \delta (\varepsilon,  x, y)  = \delta^{x}_{\varepsilon} y$ is continuous. Moreover, it can be continuously extended to $\displaystyle Cl(dom \, \delta)$ and we have 
$$\lim_{\varepsilon\rightarrow 0} \delta_{\varepsilon}^{x} y \, = \, x \quad . $$

\item[{\bf A2.}] For any  $x, \in K$, $\displaystyle \varepsilon, \mu \in \Gamma_{1}$ and $\displaystyle u \in 
\bar{B}_{d}(x,A)$   we have: 
$$ \delta_{\varepsilon}^{x} \delta_{\mu}^{x} u  = \delta_{\varepsilon \mu}^{x} u  \ .$$

\item[{\bf A3.}]  For any $x$ there is a  function $\displaystyle (u,v) \mapsto d^{x}(u,v)$, defined for any $u,v$ in the closed ball (in distance d) $\displaystyle 
\bar{B}(x,A)$, such that 
$$\lim_{\varepsilon \rightarrow 0} \quad \sup  \left\{  \mid \frac{1}{\varepsilon} d(\delta^{x}_{\varepsilon} u, \delta^{x}_{\varepsilon} v) \ - \ d^{x}(u,v) \mid \mbox{ :  } u,v \in \bar{B}_{d}(x,A)\right\} \ =  \ 0$$
uniformly with respect to $x$ in compact set. 

\end{enumerate}

\begin{rk}
The "distance" $d^{x}$ can be degenerated: there might exist  
$\displaystyle v,w \in U(x)$ such that $\displaystyle d^{x}(v,w) = 0$. 
\label{imprk}
\end{rk}

For  the following axiom to make sense we impose a technical condition on the co-domains $\displaystyle V_{\varepsilon}(x)$: for any compact set $K \subset X$ there are $R=R(K) > 0$ and 
$\displaystyle \varepsilon_{0}= \varepsilon(K) \in (0,1)$  such that  
for all $\displaystyle u,v \in \bar{B}_{d}(x,R)$ and all $\displaystyle \varepsilon \in \Gamma$, $\displaystyle  \nu(\varepsilon) \in (0,\varepsilon_{0})$,  we have 
$$\delta_{\varepsilon}^{x} v \in W_{\varepsilon^{-1}}( \delta^{x}_{\varepsilon}u) \ .$$

With this assumption the following notation makes sense:
$$\Delta^{x}_{\varepsilon}(u,v) = \delta_{\varepsilon^{-1}}^{\delta^{x}_{\varepsilon} u} \delta^{x}_{\varepsilon} v . $$
The next axiom can now be stated: 
\begin{enumerate}
\item[{\bf A4.}] We have the limit 
$$\lim_{\varepsilon \rightarrow 0}  \Delta^{x}_{\varepsilon}(u,v) =  \Delta^{x}(u, v)  $$
uniformly with respect to $x, u, v$ in compact set. 
\end{enumerate}

\begin{defi}
A triple $(X,d,\delta)$ which satisfies A0, A1, A2, A3, but $\displaystyle d^{x}$ is degenerate for some 
$x\in X$, is called degenerate dilatation structure. 

If the triple $(X,d,\delta)$ satisfies A0, A1, A2, A3, A4 and 
 $\displaystyle d^{x}$ is non-degenerate for any $x\in X$, then we call it  a 
 dilatation structure. 
 \label{defweakstrong}
\end{defi}

\subsection{Tangent bundle of a dilatation structure}
\label{induced}

The following two theorems describe the most important metric and algebraic 
properties of a dilatation structure. As presented 
here these are condensed statements, available in full length as theorems 7, 8,
10 in \cite{buligadil1}.

\begin{thm}
Let $(X,d,\delta)$ be a  dilatation structure. Then the metric space $(X,d)$ 
admits a metric tangent space at $x$, for any point $x\in X$. 
More precisely we have  the following limit: 
$$\lim_{\varepsilon \rightarrow 0} \ \frac{1}{\varepsilon} \sup \left\{  \mid d(u,v) - d^{x}(u,v) \mid \mbox{ : } d(x,u) \leq \varepsilon \ , \ d(x,v) \leq \varepsilon \right\} \ = \ 0 \ .$$
\label{thcone}
\end{thm}

\begin{thm}
Let $(X,d,\delta)$ be a dilatation structure. Then for any $x \in X$ the triple 
 $\displaystyle (U(x), \Sigma^{x}, \delta^{x}, d^{x})$ is a normed local 
 conical group. This means: 
 \begin{enumerate}
 \item[(a)]  $\displaystyle \Sigma^{x}$ is a local group operation on $U(x)$,
 with $x$ as neutral element and $\displaystyle \, inv^{x}$ as the inverse element
 function; 
  \item[(b)] the distance $\displaystyle d^{x}$ is left invariant with respect 
  to the group operation from point (a); 
 \item[(c)] For any $\varepsilon \in \Gamma$, $\nu(\varepsilon) \leq 1$, the 
 dilatation $\displaystyle \delta^{x}_{\varepsilon}$ is an automorphism with 
 respect to the group operation from point (a); 
 \item[(d)] the distance $d^{x}$ has the cone property with
respect to dilatations: foar any $u,v \in X$ such that $\displaystyle d(x,u)\leq 1$ and 
$\displaystyle d(x,v) \leq 1$  and all $\mu \in (0,A)$ we have: 
$$d^{x}(u,v) \ = \ \frac{1}{\mu} d^{x}(\delta_{\mu}^{x} u , \delta^{x}_{\mu} v) 
 \quad .$$ 
 \end{enumerate}
\label{tgene}
\end{thm}

The conical group $\displaystyle (U(x), \Sigma^{x}, \delta^{x})$ can be regarded as the tangent space 
of $(X,d, \delta)$ at $x$. Further will be denoted by: 
$\displaystyle T_{x} X =  (U(x), \Sigma^{x}, \delta^{x})$.

By using proposition 5.4 \cite{siebert} and from some topological considerations
  we deduce the following characterisation of  tangent spaces asociated to some  
  dilatation structures. The following is corollary 4.7 \cite{buligacont}. 

\begin{cor}
Let $(X,d,\delta)$ be a dilatation structure with group $\Gamma = (0,+\infty)$
and the morphism $\nu$ equal to identity. 
Then for any $x \in X$ the local group  
 $\displaystyle (U(x), \Sigma^{x})$ is locally a simply connected Lie group 
 whose Lie algebra admits a positive graduation (a Carnot group).
 \label{cortang}
\end{cor}

\subsection{Equivalent dilatation structures}

\begin{defi}
Two dilatation structures $(X, \delta , d)$ and $(X, \overline{\delta} , \overline{d})$  are equivalent  if 
\begin{enumerate}
\item[(a)] the identity  map $\displaystyle id: (X, d) \rightarrow (X, \overline{d})$ is bilipschitz and 
\item[(b)]  for any $x \in X$ there are functions $\displaystyle P^{x}, Q^{x}$ (defined for $u \in X$ sufficiently close to $x$) such that  
\begin{equation}
\lim_{\varepsilon \rightarrow 0} \frac{1}{\varepsilon} \overline{d} \left( \delta^{x}_{\varepsilon} u ,  \overline{\delta}^{x}_{\varepsilon} Q^{x} (u) \right)  = 0 , 
\label{dequiva}
\end{equation}
\begin{equation}
 \lim_{\varepsilon \rightarrow 0} \frac{1}{\varepsilon} d \left( \overline{\delta}^{x}_{\varepsilon} u ,  
 \delta^{x}_{\varepsilon} P^{x} (u) \right)  = 0 , 
\label{dequivb}
\end{equation}
uniformly with respect to $x$, $u$ in compact sets. 
\end{enumerate}
\label{dilequi}
\end{defi}

\begin{prop}
Two dilatation structures $(X, \delta , d)$ and $(X, \overline{\delta} , \overline{d})$  are equivalent  if and 
only if 
\begin{enumerate}
\item[(a)] the identity  map $\displaystyle id: (X, d) \rightarrow (X, \overline{d})$ is bilipschitz and 
\item[(b)]  for any $x \in X$ there are functions $\displaystyle P^{x}, Q^{x}$ (defined for $u \in X$ sufficiently close to $x$) such that  
\begin{equation}
\lim_{\varepsilon \rightarrow 0}  \left(\overline{\delta}^{x}_{\varepsilon}\right)^{-1}  \delta^{x}_{\varepsilon} (u) = Q^{x}(u) , 
\label{dequivap}
\end{equation}
\begin{equation}
 \lim_{\varepsilon \rightarrow 0}  \left(\delta^{x}_{\varepsilon}\right)^{-1}  \overline{\delta}^{x}_{\varepsilon} (u) = P^{x}(u) , 
\label{dequivbp}
\end{equation}
uniformly with respect to $x$, $u$ in compact sets. 
\end{enumerate}
\label{pdilequi}
\end{prop}

 The next theorem shows a link between the tangent bundles of equivalent dilatation structures. 
 
 \begin{thm} 
 Let $(X, \delta , d)$ and $(X, \overline{\delta} , \overline{d})$  be  equivalent dilatation structures. Suppose that for any $x \in X$ the distance $d^{x}$ is non degenerate. Then for any $x \in X$ and 
 any $u,v \in X$ sufficiently close to $x$ we have:
 \begin{equation}
 \overline{\Sigma}^{x}(u,v) = Q^{x} \left( \Sigma^{x} \left( P^{x}(u) , P^{x}(v) \right)\right) . 
 \label{isoequiv}
 \end{equation}
 The two tangent bundles  are therefore isomorphic in a natural sense. 
 \label{tisoequiv}
 \end{thm}

\subsection{Differentiable functions}
 
 Dilatation structures allow to define differentiable functions. The idea is to 
keep only  one relation from definition \ref{dilequi}, namely (\ref{dequiva}). 
We also renounce to uniform convergence with respect 
 to $x$ and $u$, and we replace this with uniform convergence in the "$u$" variable,  
with a conical group morphism condition for the derivative.

\begin{defi}
 Let $(N,\delta)$ and $(M,\bar{\delta})$ be two conical groups. A continuous 
 function $f:N\rightarrow M$ is a conical group morphism if $f$ is a group morphism and for any $\varepsilon>0$ and $u\in N$ we have 
 $\displaystyle f(\delta_{\varepsilon} u) = \bar{\delta}_{\varepsilon} f(u)$. 
\label{defmorph}
\end{defi}

 \begin{defi}
 Let $(X, \delta , d)$ and $(Y, \overline{\delta} , \overline{d})$ be two dilatation structures and $f:X \rightarrow Y$ be a continuous function. The function $f$ is differentiable in $x$ if there exists a 
 conical group morphism  $\displaystyle Q^{x}:T_{x}X\rightarrow T_{f(x)}Y$, defined on a neighbourhood of $x$ with values in  a neighbourhood  of $f(x)$ such that 
\begin{equation}
\lim_{\varepsilon \rightarrow 0} \sup \left\{  \frac{1}{\varepsilon} \overline{d} \left( f\left( \delta^{x}_{\varepsilon} u\right) ,  \overline{\delta}^{f(x)}_{\varepsilon} Q^{x} (u) \right) \mbox{ : } d(x,u) \leq \varepsilon \right\}   = 0 , 
\label{edefdif}
\end{equation}
The morphism $\displaystyle Q^{x}$ is called the derivative of $f$ at $x$ and will be sometimes denoted by $Df(x)$.

The function $f$ is uniformly differentiable if it is differentiable everywhere and the limit in (\ref{edefdif}) 
is uniform in $x$ in compact sets. 
\label{defdif}
\end{defi}
 
 A trivial way to obtain a differentiable function (everywhere)  is to modify the dilatation structure on the target space. 
 
 \begin{defi}
 Let  $(X, \delta , d)$ be a  dilatation structure and $f:(X, d) \rightarrow (Y, \overline{d})$ be a bilipschitz  and surjective  function. We define then the transport of $(X, \delta , d)$ by $f$, named $(Y, f*\delta , \overline{d})$, by: 
 $$\left( f*\delta\right)^{f(x)}_{\varepsilon} f(u) = f \left( \delta^{x}_{\varepsilon} u \right) . $$
 \label{ddif}
 \end{defi}

 The relation of differentiability with equivalent dilatation structures is given by the following simple proposition.

 \begin{prop}
 Let  $(X, \delta , d)$ and $(X, \overline{\delta} , \overline{d})$ be two dilatation structures and $f:(X, d) \rightarrow (X, \overline{d})$ be a bilipschitz  and surjective  function. The dilatation structures $(X, \overline{\delta} , \overline{d})$ and $(X, f*\delta , \overline{d})$ are equivalent if and only if $f$ and $\displaystyle f^{-1}$ are uniformly  differentiable. 
 \label{peqd}
 \end{prop}

 We shall prove now the chain rule for derivatives, after we elaborate a bit over the definition \ref{defdif}.

 Let $(X, \delta , d)$ and $(Y, \overline{\delta} , \overline{d})$ be two dilatation structures and $f:X \rightarrow Y$  a function differentiable in $x$. The derivative of $f$ in $x$ is a 
 conical group morphism  $\displaystyle Df(x):T_{x}X\rightarrow T_{f(x)}Y$, which means that $Df(x)$ is 
 defined on a open set around $x$ with values in a open set around $f(x)$, having the properties:
 \begin{enumerate}
 \item[(a)] for any $u,v$ sufficiently close to $x$ 
 $$Df(x)\left(\Sigma^{x}(u,v)\right) = \Sigma^{f(x)}\left(Df(x)(u), Df(x)(v)\right)  , $$
 \item[(b)] for any $u$ sufficiently close to $x$ and any $\varepsilon \in (0,1]$ 
 $$Df(x)\left(\delta^{x}_{\varepsilon} u\right) = \bar{\delta}^{f(x)}_{\varepsilon}\left(Df(x)(u)\right) , $$
 \item[(c)] the function $Df(x)$ is continuous, as uniform limit of continuous 
functions. Indeed, the relation 
 (\ref{edefdif}) is equivalent to the existence of the uniform limit 
(with respect to $u$ in compact sets)
 $$Df(x)(u) = \lim_{\varepsilon\rightarrow 0} \bar{\delta}^{f(x)}_{\varepsilon^{-1}} \left( f\left( \delta_{\varepsilon}^{x} u\right)\right) . $$ 
 \end{enumerate}
 
 From (\ref{edefdif}) alone and axioms of dilatation structures we can prove properties (b) and (c). 
 We can reformulate therefore the definition of the derivative by asking that $Df(x)$ exists as an uniform 
 limit (as in point (c) above) and that $Df(x)$ has the property (a) above. 
 
 From these considerations the chain rule for derivatives is straightforward. 
 
 \begin{prop}
 Let $(X, \delta , d)$,  $(Y, \overline{\delta} , \overline{d})$ and $(Z, \hat{\delta} , \hat{d})$ be three  dilatation structures and $f:X \rightarrow Y$  a continuous function differentiable in $x$, $g:Y \rightarrow Z$  a continuous function differentiable in $f(x)$. Then $gf:X \rightarrow Z$ is   differentiable in $x$ and 
$$D gf (x) = Dg(f(x)) Df(x) . $$
\label{pchain}
\end{prop}

\paragraph{Proof.}
Use property (b) for proving that $Dg(f(x)) Df(x)$ satisfies (\ref{edefdif}) for the function $gf$ and $x$. 
Both $Dg(f(x))$ and $Df(x)$ are conical group morphisms,  therefore  $Dg(f(x)) Df(x)$ is a conical group 
morphism too. We deduce that  $Dg(f(x)) Df(x)$ is the derivative of $gf$ in $x$. \quad $\square$

\end{document}